%

\input mssymb
\def\rest{\mathord{\restriction}}
\def\ms{\medskip}

\def\ss{\smallskip}
\def\mb{\medbreak}
\def\bb{\bigbreak}

\def\keilu{\simeq}
\def\phi{\varphi}
\def\cov{{\rm cov}}
\def\Lim{{\rm Lim\,}}
\def\model{\models}
\def\sat{\models}
\def\su{\subseteq}
\def\a{\alpha}
\def\b{\beta}
\def\c{\gamma}

\def\d{\delta}
\def\l{\lambda}
\def\k{\kappa}

\def\om{\omega}
\def\al{\aleph}

\def\lng{\langle}
\def\rng{\rangle}
\def\ov{\overline}
\def\sm{\setminus}
\def\nac{{\rm nacc }\,}
\def\nacc{\nac}
\def\acc{{\rm acc}}

\def\dom{{\rm dom}}
\def\cf{{\rm cf}}
\def\inv{{\rm Inv}}
\def\s#1{{\bf Section #1}: }

\def\epsilon{\nu}
\def\otp{{\rm otp}}
\baselineskip=16pt
\def\ran{{\rm  ran}}

\def\endproof{\hfill $\dashv$}
\def\lz{\ms}
\def\hz{\ss}
\def\new{\par\noindent}
\def\hz{\new}
\def\lz{\mb}
\def\ku{\bf}
\def\Inv{{\rm INV}}
\centerline {\bf Non-existence of Universal Orders in Many Cardinals}

\centerline {Menachem Kojman and Saharon Shelah\footnote \dag{Supported
by a BSF grant. Publication No. 409}}

\centerline{Department of Mathematics, Hebrew University}
\centerline{91904 Jerusalem, Israel}

\centerline{ABSTRACT}
{\baselineskip=12pt
\narrower\narrower
 Our theme is that not every interesting question in set
theory is independent of $ZFC$. We give an example of a first order
theory $T$ with countable $D(T)$ which cannot have a universal model
at $\aleph_1$ without CH; we prove in $ZFC$ a covering theorem from
the hypothesis of the existence of a universal model for some theory;
and we prove --- again in ZFC --- that for a large class of cardinals
there is no universal linear order (e.g. in every
$\aleph_1<\l<2^{\aleph_0}$). In fact, what we show is that if there is
a universal linear order at a regular $\l$ and its existence is not a
result of a trivial cardinal arithmetical reason, then $\l$
``resembles'' $\aleph_1$ --- a cardinal for which the consistency of
having a universal order is known. As for singular cardinals, we show
that for many singular cardinals, if they are not strong limits then
they have no universal linear order. As a result
of the non existence of a universal linear order, we show the
non-existence of universal models for all theories possessing the
strict order property (for example, ordered fields and groups, Boolean
algebras, p-adic rings and fields, partial orders, models of PA and so
on).
\lz{\bf Key words}: Universal model, Linear Order, Covering Numbers,
Club guessing, Strict Order Property
\hz 	{\bf Subject Classification}: Model Theory, Set Theory, Theory
of Orders \bigskip}
\magnification=1100

\ms
{\bf 0. Introduction} 
\mb
{\bf General Description}

This paper consists entirely of proofs in ZFC. We can even dare to
recommend reading it to anybody who is interested in linear orders or
partial orders in themselves, and to whom axiomatic set teory and model theory
are of less interest.
Such a reader should, though, consult the appendix to this paper, or a
standard textbook like [CK] for the notion of ``elementary submodel'',
and confine his reading to sections 3, 4 and 5.

The general problem addressed in this paper is the computation of the
universal spectrum of a theory (or a class of models), namely the
class of cardinals in which the theory (the class) has a universal
model. (A definition of ``universal model'' is found below). As the
universal spectrum of a theory usually depends on cardinal arithmetic,
and even on the particular universe of set theory in which a given
cardinal arithmetic holds (see below), the problem of determining the
universal spectrum of a theory must be rephrased as: under which
cardinal arithmetical assumptions {\it can} a given theory (class)
possess a universal model in a given cardinality $\l$?

All results in this paper are various {\it negative} answers in ZFC to
this question, namely theorems of the form `` if $C(\l)$ (some
cardinal arithmetic condition on a cardinal $\l$) then there is {\it
no} universal model of $T$ at cardinality $\l$''. In general, it is
harder to prove such theorems when the cardinal $\l$ in question is
singular. Such theorems are first proved for the case where $T$ is the
theory of linear orders, and then are shown to hold also for a larger
class of theories, including the theory of Boolean algebras, the
theory of ordered fields, the theory of partial orders and others.

\lz{\bf Background and detailed content} 
\new A
universal model at power $\l$, for a class of models $K$, is a model
$M\in K$ of cardinality $\l$ with the property that for all $N\in K$
such that $|N|\le \l$ there is an embedding of $N$ into $M$. At this
point let us clarify what ``embedding'' means in this paper. If
$K=MOD(T)$ is the class of models of a first order theory $T$, then
``embedding'' should be understood as ``elementary embedding'' when
$T$ is complete, and ``universal'' is with respect to elementary
embeddings; when $T$ is not complete (e.g. the theory of linear
orders, the theory of graphs or the theory of Boolean algebras),
``embedding'' is an ordinary embedding, namely a 1-1 function which
preserves all relations and operations, and ``universal'' is with
respect to ordinary embeddings. This distinction is necessary,
because there are theories for which universal models in the sense of
an ordinary embedding exist, whereas universal models in the sense of
an elementary embedding do not exist (see appendix for such an
example).

Although the notion ``universal model'' is older then its relative,
``saturated model'', and arises more often and more naturally in other
branches of mathematics, it has won less attention, perhaps because
answers to questions involving the former notion were harder to get.
As one example of a contribution to the theory of universal models we
can quote [GrSh 174], in which it is shown that the class of locally
finite groups has a universal model in any strong limit of cofinality
$\aleph_0$ above a compact cardinal. The class $\{\l: T$ has a
saturated model of cardinality $\l\}$ has been characterized for a
first order theory $T$ (See the situation --- with history --- in
[Sh-a] or [Sh-c], VIII.4).

Saturated models are universal, and their existence in known at
cardinals $\l$ such that $\l=\l^{<\l}>|T|$ or just
$\l=\l^{<\l}\ge|D(T)|$ ($D(T)$ is defined below) for every $T$;
furthermore, when $\l=2^{<\l}$, essentially the same proof gives a
``special model'', which is also universal (for these results see
[CK]).  Therefore the problem of the existence of a universal model
for a first order theory remains unsettled in classical model theory
only for cardinals $\l<2^{<\l}$.

The consistency of {\it not} having universal models at such $\l$'s for
all theories which do not have to have one at every infinite power is
very easy (see appendix). In the other direction, the second author
proved in [Sh 100] the consistency of the existence of a universal
order at $\aleph_1$ with the negation of CH, and, in [Sh 175], [Sh
175a], proved the consistency of the existence of universal graph at
$\l$, if there is a $\k$ such that $\k=\k^{<\k}<\l<2^\k=\cf(2^\k)$.
One could expect at that point to prove that every theory $T$ which
has no trivial reason for not having a universal model at
$\aleph_1<2^{\aleph_0}$, can have one.  (By a ``trivial reason'' we
mean an uncountable $D(T)$. $D(T)$ is the set of all complete
$n$-types over the empty set, $n<\omega$; it is known and easy to
prove that if $D(T)$ is uncountable, it is of size $2^{\aleph_0}$; and
every type in $D(T)$ must be realized in a universal model).  But this
is not the case.  In Section 1. we show that there is a first order
theory $T$ with $|D(T)|=\aleph_0$ (which is even
$\aleph_0$-categorical) which has a universal model in $\aleph_1$ iff
CH.

An attempt to characterize the class of theories for which it is
{\it consistent} to have a universal model at $\aleph_1<2^{\aleph_0}$ was
done by Mekler. Continuing [Sh 175]  he has shown in [M]
that it is consistent with the negation of CH that every universal
theory of relational structures with the joint embedding property and
amalgamation for $P^-(3)$-diagrams and only finitely many isomorphism
types at every finite power, has a universal model at $\aleph_1$. He
has also shown, continuing [Sh 175a], 
that it is consistent with $\k^{<\k}=\k<\cf(2^\k)<\l<2^\k$ that every
4-amalgamation class, which in every finite power has only finitely
many isomorphism types, has a universal structure in power $\l$.

In Section 2 we prove a covering theorem which shows, as one
corollary, that if $2^{\aleph_0}=\aleph_{\omega_1}$, there are no
universal models for non-$\omega$-stable theories in every regular
$\l$ below the continuum.

In Section 3 we prove in $ZFC$ several non-existence theorems for
universal linear orders in regular cardinals. We show that there can
be a universal linear order at a regular cardinal $\l$ only if
$\l=\l^{<\l}$ or if $\l=\mu^+$ and $2^{<\mu}\le\l$.  In Section 5 we
prove non-existence theorems for universal linear orders in singular
cardinals. For example, if $\mu$ is not a strong limit and is not a
fix point of the $\aleph$ function, then there is no universal linear
order in $\mu$.

In Section 5 we reduce the existence of a universal linear order in
cardinality $\l$ to the existence of a universal model for any theory
possessing the strict order property. Thus the non existence theorems
from Sections 3 and 4 which were proven for linear orders are shown to
hold for a large collection of theories.

The combined results from Sections 3, 4 and 5 show that it is
impossible to generalize [Sh 100] in the same fashion that [Sh 175a] and
[M] generalize [Sh 175]: While the proof of the consistency of having
a universal graph in $\aleph_2<2^{\aleph_0}$ generalizes the proof for
the case $\aleph_1<2^{\aleph_0}$, the 
consistency of universal linear order is true for the former case and is
{\it false} for the latter. This points out an interesting  difference
between the theory of order and the theory of graphs.

The second author is interested in the classification of unstable
theories (see [Sh 93]). With respect to the problem of determining the
stability spectrum of a theory $T$ (namely the class $K_T=\{\l:T$ has
a universal model at $\l\}$), there are several more results which
were obtained, in addition to what is published here: the main one is
a satisfactory distinction between superstabele to
stable-unsuperstable theories.
\hz{\bf Notation and Terminology}:
By ``order'' we shal mean linear order. $|M|$ denotes the universe of
a model $M$ and $||M||$ denotes its cardinality.

\bb
\s{1} {\bf A theory without universal models in
$\aleph_1<2^{\aleph_0}$}.
\ms 
We present a theory $T$. In the  language $L(T)$ there are  two
$n$-ary relation symbols, $R_n(\cdots)$ and $P_n(\cdots)$ for every
natural number $n\ge 2$. $T$ has no constants or function sumbols.
The axioms of $T$ are:
\item {1.} The sentences saying that $P_n$ and $R_n$ are invariant under
permutation of  arguments and that $P_n(x_1,\cdots,x_n)$ and
$R_n(x_1,\cdots,x_n)$ do not hold if for some $1\le i<j\le
n$, $x_i=x_j$, for all $n\ge 2$ 
\item {2.} For each $n$ the sentence saying that there are no $2n-1$
distict elements, $x_1,\cdots,x_n$, $y_1,\cdots,y_{n-1}$ such that
$P_n(x_1,\cdots,x_n)$ and for all $1\le i \le
n$, $R_n(y_1,\cdots,y_{n-1},x_i)$
\ms 
{\bf Fact 1.1}: 1. There are only finitely many quantifier free
$n$-types of $T$ for every finite $n$. \new 
2. $T$ has the joint embedding property and the amalgamtion property. 
\ss {\ku Proof}: 1. is obvious. Suppose that $M,N$ are two
models of $T$ which agree on their intersection. As $T$ is universal,
the intersection is also a model of $T$. Define a model $M'$ such that
$|M'|=|M|\cup|N|$ and such that $P_n^{M'},R_n^{M'}$ equal,
respectively, $P_n^M\cup P_n^N$ and $R_n^M\cup R_n^N$. Suppose to the
contrary that $M'$ does not satisfy $T$. So for some $n$ there are
$a,_1,\cdots,a_n,b_1,\cdots,b_{n-1}$ which realize the forbidden type.
Certainly, $\{a_1,\cdots,a_n,b_1,\cdots,b_{n-1}\}\not\su M$, as $M\sat
T$, and $\{a_1,\cdots,a_n,b_1,\cdots,b_{n-1}\}\not\su N$ as $N\sat T$.
So either
\item {(a)} there is some $a_i\notin M$ and $a_j\notin N$,
\item {(b)} there is some $a_i\notin M$ and $b_j\notin N$
or
\item {(c)} there is some $b_i\notin M$ and $b_j\notin N$.

If (a) holds this contradicts  $P(a_1,\cdots,a_n)$; if (b)
holds this contradicts $R(b_1,\cdots,b_{n-1},a_i)$; and if
(c) holds this contradicts $R(b_1,\cdots,b_{n-1},a_1)$.\endproof

\ms{\bf Fact 1.2} 
 $T$ has a universal homogeneous model $M$ at $\aleph_0$.
 \ss
This should be well known, but for completeness of presentation's sake
we give:
\ss
{\ku Proof}: Construct an increasing sequence of finite models
$M_n$:
\item{1.} $M_n\sat T$ and $M_n$ is finite.
\item{2.} $M_n\subset M_{n+1}$.
\item{3.} In $M_{n+1}\sm M_n$ all quantifier free types (of $T$) over $M_n$
are realized.

As $T$ has only finitely many quantifier free types over every finite
set, and because of fact 1.1, this construction is possible.
The model $M=\cup M_n$ clearly satisfies $T$. 

Suppose
that $h$ is a finite embedding from any other model $N$ into $M$ and
that $a\in N\sm \dom(h)$. There is some $n_0$ such that $\ran(h)\su
M_{n_0}$. In $M_{n_0+1}$ there is some $b$ such that its relational
type over $\ran(h)$ (in $M$) equals the relational type of $a$ over
$\dom(h)$ (in $N$). Set $h'=h\cup\lng a,b\rng$ to obtain an embedding
with $a$ in its domain. By this observation it is immediate that every
countable model of $T$ is embeddable into $M$. Hence the universality
of $M$ in $\aleph_0$.

 As there are no unary relation symbols in $T$, any $h=\{\lng
a_1,a_2\rng\}$, where $a_1,a_2\in M$ is an embedding. Suppose that $h$
is a finite embedding, $\dom(h),\,\ran(h)\su M$ and that $b\in M\sm
\ran(h)$. Pick, as before, some $a\in M\sm \dom(h)$ such that its
relational type over $\dom(h)$ equals the relational type of $b$ over
$\ran(h)$, and extend $h$ to include $b$ in its range. These
observations show that for every two sequences $\ov a,\ov b\in M^n$
there is an automorphism $f$ of $M$ with $f(\ov a)=\ov b$.  Hence $M$ is
homogeneous.\endproof

Denote by $T_1$ the theory $Th(M)$, the theory of the model $M$.  Then
\ms
{\bf 1.3 Fact} $T_1$ is a complete theory
extending $T$, which admits elimination
of quantifiers and is $\aleph_0$-categorical.
\ms {\ku Proof} Clearly, every simple existential formula is equivalent to a
quantifier free formula in $T_1$. Hence the elimimation of
quantifiers. By fact 1.1 there are only finitely many $n$-types 
of $T$ over the empty set. Therefore $T$ is $\aleph_0$-categorical
(see [CK] for details). \endproof
 
\ms {\bf 1.4 Fact} (1) For every infinite model $M\sat T$ there is a model $M'$
such that $M\su M'\sat T_1$ and $||M||=||M'||$. 
\new (2) $T$ and $T_1$ have the same spectrum of universal
models, namely for every cardinal $\l$, $T$ has a universal model in
$\l$ (with respect to ordinary embeddings) iff $T_1$ has a universal
model in $\l$ (with respect to elementary embeddings). 
\ss{\ku Proof}: (1) follows by the compactness theorem, as every
finite submodel of $M$ (even countable) satisfies $T$, and is
therefore embeddable into the countable model of $T_1$. For (2), we may
forget about finite cardinals, as neither of both theories has
universal finite models (in fact $T_1$ has no finite models at all). 
Suppose first that $M\sat T$ is universal for $T$ in power $\l$. Then
by 1. there is some $M'\supseteq M$, a model of $T_1$ of the same
cardinality. Let $N\sat T_1$ be arbitrary of power $\l$. As $N\sat T$,
there is an embedding $h:N\to M$. $h$ is also an embedding into $M'$.
As $T_1$ has elimination of quantifiers, $h$ is an elementary
embedding of $N$ into $M'$. So $M'$ is a universal model of $T_1$ in
$\l$. Conversely, suppose that $M\sat T_1$ is universal for $T_1$ in
power $\l$. In particular, $M\sat T$. Let $N\sat T$ be arbitrary of
power $\l$. By (1) there is some $N'\supseteq N$ of cardinality $\l$,
$N'\sat T_1$. Let $h:N'\to M$ be an elementary embedding. $h\rest N$
is an embedding of $N$ into $M$. So $M$ itself is universal for
$T$.\endproof
\new
{\bf 1.5 Remark}: 1. $T$ does not satisfy the $3$-amalgamation property, as
seen by a simple example. 
\item {2.} Also $T_1$ has the joint embedding property.

\new
{\bf 1.6 Theorem} $T$ has a universal model in $\aleph_1$ iff
$\aleph_1=2^{\aleph_0}$.
\new{\ku Proof}: If $\aleph_1=2^{\aleph_0}$ then all countable
theories have universal models in $\aleph_1$.
We proceed now to prove that $2^{\aleph_0}>\aleph_1$ implies that $T$
has no universal model at $\aleph_1$. Suppose to the contrary that CH
fails, but that $M$ is a universal model at $\aleph_1$. Without loss
of generality, $|M|=\omega_1$. We define now $2^{\aleph_0}$ models of
$T$: for each $\eta\in\,^\omega 2$ let $M_\eta$ be a model with universe
$\omega_1$ such that
\item a. $P_n^{M_\eta}=[\omega_1]^n$ iff $R_n^{M_\eta}=\emptyset$ iff
$\eta(n)=0$
\vskip0.1cm
\item b. $R_n^{M_\eta}=[\omega_1]^n$ iff $P_n^{M_\eta}=\emptyset$ iff
$\eta(n)=1$

For each $\eta$ the model $M_\eta$ trivially satisfies $T$.

As $M$ is universal, we can choose for each $\eta$ an embedding
$h_\eta:M_\eta\longrightarrow M$. Let $M^*_\eta$ be the model obtained
from $M$ by enriching it with the relations of $M_\eta$ and the
function $h_\eta$. Let $C_\eta$ be the closed unbounded set $\{\delta\in
\omega_1:M^*_\eta\rest\delta\prec M^*_\eta\}$, and let $\delta_\eta\in
C_\eta$.

As we have $2^{\aleph_0}\,\eta$'s, by the pigeon hole principle there
are more than $\aleph_1$ sequences, $\langle \eta_i:i<i(*)\rangle$,
such that for all $i<i(*)$ $\delta_{\eta_i}=\delta_0$. As there are
only $\aleph_1$ possible values to $h_{\eta_i}(\delta_0)$, we may
assume that for all $i<i(*)$, $h_{\eta_i}(\delta_0)=\gamma_0$ for some
fixed $\gamma_0$ and that $\eta_i\rest2$ is fixed. (Note that by
elementarity, and as $h$ is one to one, $\gamma_0\ge\delta_0$). Pick
now $i<j<i(*)$. There exists an $n>1$ with $\eta_i(n)\not=\eta_j(n)$,
and assume by symmetry $\eta_i(n)=0$. This means that every $n$-tuple
of distinct members of the range of $h_{\eta_i}$ satisfies the
relation $P_n^M$, while every $n$-tuple of distinct members of the
range of $h_{\eta_j}$ satisfies the relation $R_n^M$. We intend now to
derive a contradiction by constructing the forbidden type inside $M$. Pick any
$n-1$ points, $b_1,\cdots,b_{n-1}\in \delta_0$ in the range of
$h_{\eta_j}$. Notice that $M\model R_n(b_1,\cdots,b_{n-1},\gamma_0)$.
Work now in $M^*_{\eta_i}$.
\item {} $M^*_{\eta_i}\model \exists(x) R_n(b_1,\cdots,b_{n-1},h_{\eta_i}(x))$ as
$\delta_0$ witnesses this. 

So by elementarity there is such an $c_1$
below $\delta_0$ with $a_1=h_{\eta_i}(c_1)$ also below $\delta_0$.
\item {} $M^*_{\eta_i}\model \exists(x\not=c_1)
R_n(b_1,\cdots,b_{n-1},h_{\eta_i}(x))$ as $\delta_0$ witnesses this.

So by elementarity we can find $c_2,a_2$ below $\delta_0$. 
We proceed
by induction, each time picking $a_{i+1}$ differrent from all
the previous $a$'s. So when $i=n$ we have constructed the forbidden
type, as $a_1,\cdots,a_n$, being in the range of $h_{\eta_i}$,
satisfy the relation $P_n^M$. This contradicts $M\sat T$.\endproof

The proof above tells us a bit more than is stated in theorem 1.6:
what was actualy done, was to construct $2^{\aleph_0}$ models of $T$,
each of size $\aleph_1$, such that no $2^{\aleph_0}$ can be embedded
into a single model of $T$. But this construction uses no special
feature of $\aleph_1$ and the models as defined above can be defined
in any cardinality. Let us state the following.
\new{\bf 1.7 Theorem}: Let $T$ be the theory in 1.6. 
If $\aleph_0<\l=\cf\l<2^{\aleph_0}$ and
$\mu<2^{\aleph_0}$ are cardinals, then for every family of models of $T$
$\{M_i:i<\mu\}$, each $M_i$ of cardinality $\l$, there is a model $M$
of $T$ which cannot be embedded into any $M_i$ in the family.
\new{Proof}: Suppose such a family is given.  
As in the previous proof, there are $2^{\aleph_0}$ trivial models
of $T$, each with universe $\l$. Suppose that each of them is embedded
into some member of the family. As $\mu<2^{\aleph_0}$, there must be
a fixed member $M_{i(*)}$ of the family into which more than $\l$ such
models are embedded. The contradiction follows now as above.\endproof

\bb
\s{2} {\bf A covering theorem}

We prove here a theorem that as one consequence puts a restriction on
the cofinality of $2^\l$ --- provided there is a universal model for a
suitable theory in some cardinality $\k\in (\l,2^\l)$.
\ms 
{\bf 2.1 Theorem } Let $T$ be a first order theory,
$\lambda<\kappa<\mu$ cardinals. Suppose that $T$ has a universal model
at $\kappa$ and that there is a model $M$ of $T$, $|M|=\mu$, with a
subset $A\subseteq
\mu$, $|A|=\l$ such that $|S(A)|\ge\mu$, namely there are $\mu$
complete 1-types
over $A$.  {\it Then}
there is a family $\langle B_i:i<\kappa^\l\rangle\subseteq
[\mu]^\kappa$ which covers $[\mu]^\kappa$, namely for every $C\in
[\mu]^\kappa$ there is an $i<\kappa$ such that $C\subseteq B_i$.

As corollaries we get
\ms{\bf 2.2 Corollary }: If $cf(2^\lambda)\le \kappa <2^\lambda$ and $T$ is
a first order theory  possessing the independence property then $T$
has no
universal model at $\kappa$.
\lz
{\bf 2.3 Corollary }: Suppose that $2^{\aleph_0}=\aleph_{\omega_1}$. Then  all
theories unstable in $\aleph_0$ (e.g. the theory of graphs, the theory of
linear order and so on) do not have a universal model in any cardinal
$\kappa\in [\aleph_1,\aleph_{\omega_1})$.
\ms
{\ku Proof of Corollary 2.2 from Theorem 2.1}: If $T$ possesses the
independence property, then there is a model of $T$ in which
$2^\lambda$ types over a set of size $\lambda$ are realized. By
Theorem 2.1, if there were a universal model for  $T$ at
$\kappa$, then there would be a covering family of
$[2^\lambda]^\kappa$ of size $2^\lambda$; but as $cf(2^\lambda)\le
\kappa$, this is clearly impossible.\endproof
\ss
{\ku Proof of Corollary 2.3 from Theorem 2.1}: If a countable theory
$T$ is not stable at $\aleph_0$ then $T$ has a model in which
$2^{\aleph_0}$ complete 1-types over a countable set are realized. So
a universal model at $\kappa\in[\omega_1,2^{\aleph_0})$ would imply,
by Theorem 2.1, that there is a covering family of
$[2^{\aleph_0}]^\kappa$ of size $2^{\aleph_0}$, which is impossible as
$\cf(2^{\aleph_0})\le\kappa$.\endproof
\ms
{\ku Proof of theorem 2.1}: Let $U$ be a universal model for $T$ with
universe $\kappa$ and let $M$ be a model of $T$ with a subset
$A\subset |M|$ of size $\lambda$ with $\langle p_i\in
S(M):i<\mu\rangle$ a sequence of distinct complete types over $A$.
Without loss of generality $M$ is of size $\mu$ and in $M$ all $p_i$'s
are realized. We can further assume (by enumerating $|M|$) that
$|M|=\mu+\mu$, that $p_i$ is realized by the element $i$ and that
$A=\{\alpha:\mu\le\alpha<\mu+\lambda\}$. For each submodel $N$, $A\cup
B\subseteq N\prec M$ of size $\kappa$  pick an embedding
$h_N:N\to U$. For each function $f:A\to \kappa$
let $C_f$ be the set of submodels $\{N\su M:|N|\le \kappa{\rm, }\,A\su N
{\rm, and }\; h_N\rest A=f\}$.
\hz
{\bf 2.4 Claim}: For each $f\in \kappa^A$, $|\cup C_f\cap \mu|\le \kappa$.
\ss
{\ku Proof}: Enumerate all members of $C_f$ in a sequence
$\langle N_\alpha:\alpha<\alpha(*)\rangle$ and define a function $g:\cup
C_f\cap\mu\to \kappa$ by induction on $\alpha$ as follows: 
$g\rest
(N_\alpha\setminus\cup_{\beta<\alpha}N_\b)=h_{N_\alpha}\rest
(N_\alpha\setminus\cup_{\beta<\alpha}N_\b)$. We are done if $g$ is a 1-1
function. This follows from
\lz{\bf 2.5 Fact} If $i\in N_{\alpha},\,j\in
N_{\beta}$, and  $i<j<\mu$
then $h_{N_\alpha}(i)\not=h_{N\beta}(j)$.
\lz {\ku Proof of Fact}: As both embeddings aggree on the image of
$A$ and $i,j$ realize different types over $A$, the fact is immediate.
\hz
{\bf 2.6 Claim}: The family $\{(\cup C_f)\cap\mu:f\in\,^A\kappa\}$ is a
covering family of $[\mu]^\kappa$ of size $\lambda^\kappa$.
\lz{\ku Proof of claim}: Clearly the size of the family is as stated.
Let $B\in[\mu]^\kappa$ be any set. Then it is a subset of some
elementary submodel  $N\subseteq M$ which contains $A$ as a subset. So it is a subset of
$\cup C_{f_N\rest A}\cap \mu$.\endproof  
\bb\s{3} {\bf Non-existence of universal linear orders}

In this section we prove some non-existence theorem for universal
linear orders in regular cardinals. We start by showing that there is
no universal linear order in a regular cardinal $\l$ if
$\aleph_1<\l<2^{\aleph_0}$.  We shall generalize this for more regular
cardinals later in this section. The combinatorial tool which enables
these theorems is the guessing of clubs which was introduced in [Sh-e]
and can be found also in [Sh-g], which will, presumably, be available
sooner. Proofs of the relevant combinatorial principles are repeated
in the appendix to this paper for the reader's convenience.
\hz{\bf 3.1 Definitions}
\item {1.} If  $C$ is a set of ordinals,  $\d$ an ordinal,  we denote by
$\d_C^s$ the element ${\rm min}\{C\setminus (\d+1)\}$ when it exists. 
\item {2.} a {\it cut} $D$ of a linear order $O$ is pair $\lng
D_1,D_2\rng$ such that $D_1$ is an initial segment of $O$, namely
$D_1\su |O|$ and $y<x\in D\Rightarrow y\in D_1$, $D_2$ is an end
segment, namely $D_2\su |O|$ and $y>x\in D_2\Rightarrow y\in D_2$, $D_1\cap
D_2=\emptyset$ and $D_1\cup D_2=|O|$. If $O_1\su O_2$ are linear
orders, then an element $x\in O_2\sm O_1$ {\it realizes} a cut $D$ of
$O_1$ if $D_1=\{y\in O_1:y<x\}$.
\item  {3.} Let $O=\cup_{j<\l}O_j$ be an increasing continuous union of
linear orders,  let
$\d\in \l$ be limit, and let $C\su\d$ be unbounded in $\d$. Let $x\in
(O\sm \cup_{j<\d}O_j)$. Define $\inv_O(C,\d,x)$, {\it the invariance of $x$
in $O$ with respect to} $C$, as $\{\a\in C:\exists y\in
O_{\a_C^s}$ such that $y$ and $x$ realize the same cut of
$O_\a\}$. Note that this definition is applicable also to cuts
(rather than only to elements).
\item {4.} A {\it $\k$-scale} for $\l$ is a sequence
 $\ov C=\lng c_\d:\d\in S\rng$
where $S=\{\a<\l:\cf\a=\cf\k\}$ and for every $\d\in S$, $c_\d$ is a club
of $\d$ or order type $\k$, $c_\d=\lng \a^\d_j:j<\k\rng$ is an
increasing enumeration of $c_\d$. If
$O=\cup_{i<\l} O_i$ is a linear order represented as a continuous
increasing union of smaller orders, and $\ov C$ is a $\k$-scale for
some $\k<\l$, let $\Inv(O,\ov C)=^{def}\{X\su \k:(\exists \d\in
S)(\exists x>\d)\inv(c_\d,\d,x)=\{\a^\d_j:\a\in X\}\}$. So $\Inv(O,\ov C)$ is
the set of all subsets of $\k$ which are obtained as an invariance of
some element in $O$ with respect to some $c_\d$ in the scale.

\hz{\bf 3.2 claim}: Suppose $h:O_1\rightarrow O_2$ is an embedding of
linear orders, $||O_1||=||O_2||=\l=\cf\l>\aleph_0$. 
Then for any representations
$O_i=\cup_{j<\l}O^i_j$, $i=1,2$, the union increasing and continuous
and each $|O^i_j|<\l$, there exists a club $E\su \l$ such that for
any $\d<\l$ and $C\su\d$ a club of $\d$ which satisfies $C\su E$,
we have 
	$$ (*)\quad (\forall x\in O_1\sm O_\d)
(\inv_{O_1}(C,\d,x)=\inv_{O_2}(C,\d,h(x))$$	
\hz
{\ku Proof of Claim}: without loss of generality we may assume that
$|O_1|=|O_2|=\l$. Define the model $M=\lng \l,<_{O_1},<_{O_2},\in,
h\rng$. Let $E=\{ \d<\l:M\rest\d\prec M$ and for all $\d\in E$,
$\d=\cup_{j<\d}O^i_j$, $i=1,2\}$. Let $x\in (O_1\sm O^1_\d)$. Note
that by elementarity $h(x)\in (\l\sm\d)$.  Suppose first that $\a$
belongs to the left hand side of the equality in $(*)$ and let $y\in
[\a,\a_C^s)$ demonstrate this. So $x$ and $y$ realize the same cut of
$O_1\rest\a$.  As $h$ is an embedding, $h(x),h(y)$ satisfy the same
cut of $h''(O_1\rest\a)$ (which equals, by elementarity,
$(h''(O_1))\rest\a$).  If $h(x),h(y)$ satisfy also the same cut of
$O_2\rest \a$ we are done, but the problem is, of course, that $h$ is
not necessarily onto.  Otherwise suppose that (w.l.o.g) $h(y)<h(x)$
and that there is an element $z\in O_2\rest\a$ such that
$h(y)<_{O_2}z<_{O_2}h(x)$.  Define in $M$ the set
$D=\{t:$  there is no $ q$ such that
$z\le_{O_2}h(q)\le_{O_2} t\}$. $D$ is definable in $M$ with parameters
in $M\rest \a^s$. By elementarity the definition is absolute between
$M$ and $M\rest \a^s$, that is $D\cap \a^s$ is the same as $D$
interpreted in $M\rest\a^s$. $D$ is
a cut of $O_2\rest \a$.  Let $D'$ be $D\cap\a$. $D'$ is definable in $M\rest\a^s$.
\hz {\bf 3.2.1 subclaim}: $h(x)$
satisfies the cut $D'$ determined by $D$.
\hz{\ku Proof of Subclaim}: let $z<_{O_2}\b<_{O_2}h(x)$. 
As there are no points in
the range of $h\rest \a$ between $z$ and $h(x)$, there are certainly none
between $z$ and $\b$. So $\b\in D'$. Conversely, suppose
$h(x)<_{O_2}\b$. Then $M$ satisfies that there is an image under $h$
(namely $h(x)$) between $z$ and $\b$. By elementarity there is such an
image $h(x')$ where $x'\in\a$. So $\b\notin D'$.
\hz As $D'$ is a cut of $O'\rest\a$ definable in $M\rest \a^s$ which is realized by
$h(x)$, elementarity assures us that is it is realized by some
$y'\in[\a,\a_C^s)$. So $\a$ belongs to the right hand side of $(*)$.
\hz
Assume that $\a$ belongs to the right hand side of $(*)$.
Then there is an element $y\in [\a,\a^s)$ which satisfies the same cut
of $O_2\rest \a$ as $h(x)$. If $y=h(y')$ for some $y'$ we are done.
Else, we note that the cut of $O_2\rest\a$ which $y$ determines is
definable in $M\rest\a^s$. Now clearly $h(x)$ and $y$ satisfy the same
cut of $O_2\rest\a$. By elementarity there is an element $y'$ such that
$h(y')$ satisfies the same cut as $y$, therefore as $h(x)$. In other
words, $\a$ belongs to the left hand side of $(*)$.\endproof

\new
{\bf
3.3 Fact} If $O$ is an order with universe $\l$ and $\ov C$ is a
$\k$-scale, then $|\Inv(O,\ov C)|\le \l$
\new
{\ku Proof}: Trivial.
\new
{\bf 3.4 Lemma} (the construction lemma): If $\l<2^{\aleph_0}$ is a regular
uncountable cardinal, $\ov C$ is an $\om$-scale, and $A\su \om$ is
given, then there is an order $O$ with universe $\l$,
$O=\cup_{i<\l}O_i$, increasing continuous union of smaller orders,
such that for every $\d<\l$ with $\cf\d=\aleph_0$,
$\inv(c_\d,\d,\d) =\{\a^\d_n:n\in A\}$.
\new{\ku Proof}:
 We define by induction on $0<\a<\l$ an order $O_\a$ with the
properties listed below. We denote by $Q$ the order of the rationals.
If $O_1\su O_2$ are linear orders, $D$ a cut of $O_1$ and $D'$ a cut
of $O_2$, we say that $D'$ {\it extends} $D$ if $D_1'\rest |O_1|=D_1$
and $D_2'=D_2$. Also note that if $O_1\su O_2$ are linear orders and
$D^1$ is a cut of $O_1$ which is not realized in $O_2$ then it
corresponds naturally to a cut $D^2$ of $O_2$. In such a case we say
that $D^1$ is (really) also a cut of $O_2$.
\item {1.}  $O_\a$ has universe $|O_\a|\in \l$.
\item {2.} If $\b+1=\a$ and $x\in (O_\a\setminus O_\b)$, then $\{y\in
O_\a\sm O_\b:x$ and $y$ satisfy the same  cut of $O_\b\}$ has
order type $Q$.
\item {3.} If $\a<\b<\c$, $\c$ is a successor, and there is a cut $D$ of $O_\a$ which is
realized by an element of $O_\b$ but is not realized by no element of
$O_\epsilon$ for $\a<\epsilon<\b$, then there is a cut $D'$ of $O_\b$,
which extends $D$, which is realized in $O_\c$ but is not realized in
$O_\epsilon$ for all $\b<\epsilon<\c$. Also for every successor $\a$
there is a cut of $O_0$ which is realized in $O_\a$ but is not realize
in $O_\b$ for every $\b<\a$.
\item {4.} If $\a$ is limit then $O_\a=\cup_{\b<\a}O_\b$.
\item {5.} If $\cf\d=\aleph_0$ and for all $\b\in
C_\a$, $|O_\b|=\b$, {\it then}
$\inv_{O_{\d+1}}(c_\d,\d,\d)=\{\a^\d_n:n\in A\}$.

There should be no problem taking care of 1--4. Assume that the
conditions of 5. are satisfied. We wish to define the order
$O_{\a+1}$.   Let
$C_\a=\lng \b_n:n<\omega\rng$. By induction on $A=\lng a_n:n<\om\rng$
define an increasing sequence of cuts, $\lng D_{\a_n}:n\in \om\rng$ such
that $D_{a_n}$ is a cut of $O_{a_n}$ which is realized for the first
time in $O_{a_{n+1}}$. Demand 3. enables this. In $O_{\a+1}$ let
$\a$ satisfy $\bigcup D_n$ to get 5.\endproof

We are almost ready to prove the non existence of a universal order in
a regular $\l$, $\aleph_1<\l<2^{\aleph_0}$.  We recall from [Sh-e]
chapter III.7.8 (see also [Sh-g]),
\hz
{\bf 3.5 Fact}: If $\l>\aleph_1$ is regular, then there is a sequence
$\ov C=\lng c_\d:\d<\l,\;\cf\d=\aleph_0\rng$, such that $c_\d\su\d$ is
a club of $\d$ of order type $\om_0$, with the property that for every
club $E\su \mu$ the set $S_E=\{\d<\l:\cf\d=\aleph_0\,{\rm and}\,c_\d\su
C\}$ is stationary.\endproof

A proof of this fact is found in the appendix.
\new {\bf 3.6 Theorem} If $\aleph_1<\l=\cf\l<2^{\aleph_0}$, then there
is no universal order in cardinality $\l$.
\new {\ku Proof}: Suppose to the contrary that $UO$ is a universal
order in cardinality $\l$. Without loss of generality, $|UO|=\l$. Fix
some club guessing sequence $\ov C=\lng
c_\d:\d<\l,\,\cf\d=\aleph_0\rng$. This is known to exist by the
previous fact. As $|\Inv(UO,\ov C)|\le\l$, there is some $A\su \om$,
$A\notin \Inv(UO,\ov C)$. Use the construction lemma to get an order
$M$ with universe $\l$ and with the property that for every $\d<\l$,
$\cf\d=\aleph_0$ implies that $\inv_M(c_\d,\d,\d)=\{\a^\d_n:n\in A\}$.
Let $h:M\to UO$ be an embedding.  Let $E_h$ be the club given by 3.2.
As $\ov C$ guesses clubs, there is some $\d(*)$ with $c_{\d(*)}\su
E_h$. Therefore, by 3.2,
$\inv_M(c_{\d(*)},\d(*),\d(*))=\inv_{UO}(c_{\d(*)},\d(*),h(\d(*)))$.
But $\inv_M(\c_{\d(*)},\d(*),\d(*))=\{\a^\d_n:n\in A\}$. This means
that $A\in \Inv(UO,\ov C)$, a contradiction to the choice of $A\notin
\Inv(UO,\ov C)$.\endproof

We wish now to generalize Theorem 3.6 by replacing $\om_0$ by a more
general $\k$. As the proof of 3.6 made use of both club guessing and
the construction lemma, we should see what remains
true  of these two facts
for $\k>\aleph_0$. The proof of the construction lemma does not
work when replacing $\aleph_0$ by some other cardinal.  We need some
extra machinery to handle the limit points below $\k$.
\mb{\bf 3.7 Lemma} (the second construction lemma) Suppose
$\k<\l=\cf\l$ are cardinals, $2^\k\ge \l $ and that there is a
stationary $S\su\l$ and sequences $\lng c_\d : \d\in S\rng $ and $\lng
P_\a : \a<\l\rng$ which satisfy:
\item {(1)} $\otp c_\d =\k$ and $\sup c_\d=\d$;
\item {(2)} $P_\a\su {\cal P}(\a)$ and $|P_\a|<\l$;
\item {(3)} if $\a\in \nacc c_\d$ then $c_\d\cap\a\in\cup_{\b<\a}P_\b$,
\ss
THEN when given such sequences and a closed $A\su {\rm Lim}\k$ there is a
linear order $O$ with universe $\l$ with the property that for every
$\d\in S$, $\inv(c_\d,\d,\d)=A_\d$, where $A_\d$ is the subset of
$c_\d$ which is isomorphic to $A$.
\ss
{\ku Proof} We pick some linear order $L$ of cardinality smaller than
$\l$ which has at least $\l$ cuts. We assumme, without loss of
generality, that $P_\a\su P_\b$ whenever $\a<\b$, that for limit $\a$
$P_\a=\cup_{\b<\a}P_\b$ and that if $\a\in
\nacc c_\d$ then $A_\d\cap \a\in P_\a$. Next we construct by induction
on $\a<\l$ an order $O_\a$ and a partial function $F$ with the following
demands:
\item {(1)} the universe of $O_\a$ is an ordinal below $\l$.
\item {(2)} $\a<\b\Rightarrow O_\a\su O_\b$, and if $\a$ is limit, then
$O_\a=\cup_{\b<\a}O_\b$. 
\item {(3)} If $x\in O_\b\sm O_\a$, then the order type of $\{y\in
O_\b: x$ and $y$ satisfy the same cut of $O_\a\}$ contains $L$ as a
suborder. Also, if $\a$ is a successor, then there is an element in
$O_\a$ which satisfies a new cut of $O_0$.
\item {(4)} If $\a<\b<\gamma$ and $\gamma$ is a succsessor, then if
$D$ is a cut of $O_\a$ which is realized in $O_\b$ but not in an
earliert stage, then there is a cut $D'$ of $O_\b$ which extends $D$
and is realized in $O_\gamma$ but is not realized in $O_\d$ for any
$\d<\gamma$.
\item {(5)} $F$ is a partial function, $\dom F\su S\times(\l\sm\Lim
\l)$. A pair $\lng\d,\a\rng\in\dom F$ iff $\a<\d$,
$\emptyset\not=A_\d\cap\a\in P_\a\sm P_{\a-1}$.  $F(\d,\a)$ is a pair
$\lng \b(\d,\a), D(\d,\a)\rng$, where $\b<\a$ and $D$ is a cut of
$O_\b$ which is realized in $O_\a$. If $\b$ is not a limit of $A_\d$
then $D$ is not realized in $O_\gamma$ for any $\gamma<\a$.
$F(\d,\a)$ depends only on $A_\d\cap \a$, namely if
$A_{\d_1}\cap\a=A_{\d_2}\cap\a$ then $F(\d_1,\a)=F(\d_2,\a)$. If
$\a<\gamma$ and $F(\d,\a),F(\d,\gamma)$ are both defined, then
$\b(\d,\a)<\b(\d,\gamma)$ and $D(\d,\gamma)$ extends $D(\d,\a)$.
\item {(6)} If $\d\in S$ then $\inv_{O_{\d+1}}(c_\d,\d,\d)=A_\d$.

As $O_0$ we pick $L$. When $\a$ is limit, we define $O_\a$ as the
union of previous orders. When $\a$ is a successor we add less then
$\l$ elements to take care of demands 3 and 4. If $A_\d\cap\a\in
P_\a\sm{\a-1}$ we must define $F(\d,\a)$. If $A_\d\cap\a$ contains
exactly one member, let $\b(\d,\a)=0$ and as $D(\d,\a)$ pick (by(3)) a
cut of $O_0$ which is realized in $O_\a$ but is not realize in
$O_\gamma$ for any $\gamma<\a$. In case the order type of
$\{\gamma<\a: F(\d,\gamma)$ is defined $\}$ is limit, we let
$D(\d,\a)=\cup_{\gamma<\a}D(\d,\gamma)$ and
$\b(\d,\a)=\cup_{\gamma<\a}\b(\d,\gamma)$. Note that $\b(\d,\a)<\a$,
because it is limit. Add more elements to $O_\a$ to realize $D$. Since
$|P_\a|<\l$, this requirement of addition of elements is satisfied by
adding less than $\l$ new elements. In
case there is a last $\gamma<\a$ for which $F(\d,\gamma)$ is defined,
let this $\gamma$ be $\b(\d,\a)$ and pick (by (4)) a cut $D$ or
$O_\gamma$ which extends $D(\d,\gamma)$ and is realized in $O_\a$, but
is not realised earlier, as $D(\d,\a)$. When $\a=\d+1$ and $\d\in S$,
let the element $\d$ realize, in $O_\a$ the cut
$\cup_{\gamma<\a}F(\d,\gamma)$. 

Having added less than $\l$ new elements, we fulfill demand (1). (2)
and (3) are obvious, and (4) and (5) have been taken care of.
\ss{\bf Claim}: demand (6) holds.
\ss{\ku Proof}: Suppose that $\d\in S$. We show by induction that for
every $x\in A_\d$, for every $y\le x$ in $c_\d$ there is some
$\gamma<y^s_{c_\d}$ which satisfies the cut of $\d$ over $O_y$ iff
$y\in A_\d$. Suppose $x$ is the first member of $A_\d$. Then the first
$\gamma$ for which $F(\d,\gamma)$ is defined satisfies
$x<\gamma<x^s_{c_d}$ by the assumptions on $\lng P_\a:\a<\l\rng$.
 $F(\d,\gamma)$ is  a cut of $O_0$
which is realized in $O_\gamma$ but not before. If $y\in c_\d\cap
x$, as $x$ is a limit of $c_\d$, $y^s_{c_d}<x$. The cut of $\d$ over
$O_0$ is $F(\d,\gamma)$. And the cut of $\d$ over $O_y$ extends this
cut. As $F(\d,\a)$ is not realized by the stage $O_y$, certainly the
cut of $\d$ over $O_y$ is not realized by this stage either. So
$y\notin \inv(C_\d,\d,\d)$. As the cut of $\d$ over $O_0$ is not
realized in $O_x$, it is really also a cut of $O_x$. This cut is
realized in $O_\gamma$, where $\gamma<x^s_{\c_d}$. So by definition,
$x\in \inv (c_\d,\d,\d)$.

In the case $x$ is a successor of $A_\d$, denote by $z$ its
predecessor in $A_\d$. The minimal $\gamma$ above $x$ for which
$F(\d,\gamma)$ is defined is smaller than $X^s_{\c_\d}$, and
$\b(\d,\gamma)$ is in the interval $(z,z^s_{c_d})$. The same argument
as in the previous case shows that for every $y\in (z,x]$, $y\in
\inv(\c_\d,\d,\d)$ iff $y\in A_\d$. When $x$ is a limit of $A_\d$, by
the induction hypothesis, for every $y<x$ the required holds. As for
$x$ itself, if $\gamma$ is the minimal above $x$ for which
$F(\d,\gamma)$ is defined, $\gamma<x^s_{c_\d}$ and $F(\d,\gamma)$ is
realized in $O_\gamma$. Therefore $x\in\inv(c_\d,\d,\d)$.\endproof

By [Sh 420] we know:
\ss{\bf 3.8 Fact} If $\k$ is a cardinal and  $\k^+<\l=\cf\l$, then
there is a stationary set $S$ and sequences $\lng c_\d:\d\in
S\rng,\quad\lng P_\a:\a\in\l\rng$ as in the assumptions of 3.7.

What is still lacking is the appropriate club guessing fact, which we
quote now from [Sh~g]:
\ms{\bf 3.9  Fact}: If $\k$ is a cardinal, $\k^+<\l=\cf\l$ and there is a
stationary set $S\su\l$ and sequences $\lng c_\d:\d\in S\rng$, $\lng
P_\a:\a<\l\rng$ as in 3.7, then there are such with the additional
property that $\lng c_\d:\d\in S\rng$ guesses clubs.\endproof
\mb {\bf 3.10 Theorem} Suppose $\l=\cf\l$ and there is some cardinal $\k$
such that $\k^+<\l<2^\k$, then there is no universal linear order in
cardinality $\l$.\endproof
\ss
{\ku Proof}: Suppose $O$ is any order of cardinality $\l$, and assume
without loss of generality that its universe is $\l$. Pick a
stationary set $S$ and sequences as in 3.7, with the a property that
$\ov C=\lng c_\d:\d\in S\rng$ guesses clubs. Pick a closed set $A\su
\Lim \k$ which is not in $\Inv(O,\ov C)$. Use 3.7 to construct an
order $O'$ with unierse $\l$ and the property that for every $\d\in
S$, $\inv_{O'}(c_\d,\d,\d)\keilu A$. If $O'$ where embedded into $O$,
some $c_\d$ would guess the club of the embedding, what would lead to a
contradiction.  So $O'$ is not embeddable into $O$, and therefore
there is no universal linear order in $\l$.
\endproof

\bb
\s{4} {\bf Singular cardinals}

We shall state now a theorem which concerns the non-existence of
universal linear orders in singular cardinals. Let us note, though,
the following well known fact first:
\hz
{\bf 4.1 Fact} If $\mu$ is a strong limit, then for every first order
theory $T$ such that $|T|<\mu$ there is a special model of size $\mu$,
and therefore also a universal model in $\mu$. \endproof

For the defintion of special model see the appendix. A special model
is universal. For more details see [CK] p. 217.

This means that for non existence of universal models we must look at
singulars which are not strong limits. We will see at the end of this
section that if, e.g., $\aleph_\omega$ is not a strong limit, then
there is no universal linear order at $\aleph_\omega$.

We recall from [Sh-g 355,5]
\hz{\bf 4.2 Definition} $\cov(\l,\mu,\theta,\sigma)$ is the minimal size of a
family $A\su[\l]^{<\mu}$ which satisfies that for all
$X\in[\l]^{<\theta}$ there are less than $\sigma$ members of $A$ whose
union covers $X$.

\hz{\bf 4.3 Theorem}: Suppose $\theta=\cf\theta<\theta^+<\k$ are regular
cardinals, $\k<\mu$ and there is a binary tree $T\su {}^{<\theta}2$ of
size $<\k$ with $>\mu^*:= \cov(\mu,\k^+,\k^+,\k)$ branches of length
$\theta$.  Then 
\new{ $(*)_{\mu,\k}$} There is no linear order of size $\mu$
which is universal for linear orders of size $\k$ (namely that every
linear order of size $\k$ is embedded in it).
\hz
{\ku Proof}: Let $\ov A=\lng A_i:i<\mu^*\rng\su [\mu]^{<\k^+}$
demonstrate the definition of $\mu^*$. Without loss of generality,
$|A_i|=\k$ for all $i$. Suppose to the contrary that there is an order
$UO=\lng \mu,<_{UO}\rng$ into which every order of size $\k$ is embedded.
Let $M_i$ be $UO\rest A_i$ for every $i<\mu^*$.  Then every $M_i$ is
isomorphic to some $M'_i$ with universe $\k$, and for every order $O$ of
size $\k$ there is a set $J\su \mu^*$, $|J|<\k$ such that $O$ is
embedded into $\cup_{i\in J}M_i$.

We fix a club guessing sequence, $\ov C=\lng c_\d:\d<\mu\quad {\rm  and }
\quad\cf\d=\theta\rng$, and an increasing continuous sequence 
$\lng P_\a:\a<\k\rng$ such that 
$P_\a$ is a family of  subsets
of 
$\a$, $|P_\a|<\k$ and for all $\a<\mu$ if $\d\in S$ and $\a\in\nac
c_\d$, 
then $(c_\d\cap\a)\in P_\a$. 
For the existence of these, see [Sh 420].

For each $\d\in S$ enumerate $c_\d$ as $\lng a^\d_i:i<\theta\rng$ in
an increasing continuous fashion. Now $T$ can be viewed as $T_\d$, a
tree of subsets of $c_\d$.  Under the assumptions we already have, it
is no lose of generality to assume that for every $\a\in\k$, if
$\a\in\nac (c_\d)$, then $T_\d\cap {\cal P}(\a)\su P_\a$. The reason is
that there are $\theta$ possibilities for the unique $i$ such that
$\a=\a^\d_i$, and for each such possibility there are $<\k$ subsets in
$T\cap {\cal P}(i)$; So we can add all the required sets into $P_\a$
without changing the fact that $|P_\a|<\k$.

So by now we have the assumptions of 3.7. Using it we construct a
linear order $O$ on $\k$, with $A\su \k$ not in 
$\{\inv_{M'_i}(c_\d,\d,x):i<\mu^*,x\in \k\}$.

Suppose now that there is an embedding $h:O\rightarrow UO$. The 
image of $h$ is covered by $\cup\{O_i:j\in J\}$ for some $J$ of size
$<\k$. Let $S_j=\{x\in \k:h(x)\in M_j\}$. Then there is some $j_0$ such
that $S_{j_0}\notin \ov {id}^a(\ov C)$, (the latter is the ideal of
non-guessing, namely $X\in \ov {id}^a(\ov C)$ iff there is a club $E$
such that $\forall(\d\in S\cap X)(c_\d\not\su E)$. This ideal is
clearly $\k$-complete.)

Let $O'$ be $O\rest S_{j_0}$. Let $O'_i=O'\rest i$ for $i<\k$. Then
this is a presentation of $L'$ as an increasing continuous union of
small orders. By 3.2, and the fact that the identity map embeds $O'$
into $O$, almost everywhere the invariance with relation to $O$ is the
same as with relation to $O'$. So we can get again the same
contradiction as in previous proofs by inspecting the embedding $h\rest
O'$.\endproof

We wish now to obtain the same results using more concrete
assumptions. We first review some facts concerning covering numbers.

Recall the well known (see e.g. [Sh-g 355,5])
\hz{\bf 4.4 Fact} If $\d<\k=\cf\k<\mu=\aleph_\d$ then $\cov(\mu,\k^+,\k^+,\k)=\mu$.
\new
{\ku Proof}: By induction on $\chi$, a cardinal, $\k<\chi\le\mu$.
\item{(a)} $\chi=\theta^+$.

For every $\alpha<\chi$ fix a family $P_\alpha\su[\alpha]^\k$ with the
proprty that for every set $A\in [\alpha]^\k$ there is a set $X\su
P_\alpha$, $|X|<\k$ and $A\su \cup X$. Let $P$ be the union of
$P_\alpha$ for $\alpha<\chi$. The size of $P$ is clearly $\chi$, and
clearly for every set $A\su\chi$n of size $\k$ there is a covering of
$A$ by less than $\k$ members of $P$.

\item{(b)} $\chi=\aleph_\b$ is a limit cardinal.

As $\mu=\aleph_\d$ with $\d<\k$, certainbly $\b<\k$. Let
$\lng\chi_i:i<\cf\b\rng$ be increasing and unbounded below $\chi$. Let
$P_i$ demonstrate that $\cov(\chi_i,\k^+,\k^+,\k)=\chi_i$, and let
$P=\cup_iP_i$. Then $|P|=\chi$. If $A\in[\chi]^\k$, cover
$A\cap\chi_i$ by less than $\k$ members of $P$. Thus to cover $A$ we
need less than $\k$ members of $P$.
\hz
{\bf 4.5 Improved fact}

If $\mu$ is a fix point of the first order (i.e. $\l=\aleph_\l$), but
not of the second order, i.e. $|\{\l<\mu:\l=\aleph_\l\}|=\sigma<\mu$,
and $\sigma+\cf\mu<\k<\mu$, then $\cov(\mu,\k^+,\k^+,\k)=\mu$.
\ss
{\bf Proof}:
Suppose that $\k<\chi<\mu$ and $\cf\chi=\k$. By the assumptions,
$\chi\not=\aleph_\chi$, say $\chi=\aleph_\d$. By [Sh 400], section 2,
${\rm pp}(\chi)<\aleph_{|\d|^{+4}}<\mu$. By [Sh-g 355, 5.4], and the
fact that $\chi$ is arbitrarily large below $\mu$, we are done.\endproof

We see that we can have arbitrarily large $\k$ below a singular $\mu$
with $\mu=\cov(\mu,\k^+,\k^+,\k)$, when $\mu$ is a limit which is not
a second order fix point of the $\aleph$ function.  But for applying
theorem 4.3 we need also a binary tree of height and size $<\mu$ with $>\mu$
branches. This happens if there is some $\sigma<\mu$ with
$2^{<\sigma}<\mu$ and $2^\sigma>\mu$.  So we can state

\hz{\bf 4.6 Corollary}: If $\mu$ is a singular cardinal which is not a second
order fix point, and there is some $\sigma<\mu$ such that
$2^{<\sigma}<\mu<2^\sigma$, then there is no universal linear order of
power $\mu$.
\new
{\ku Proof}: Let $T$ be the tree $^{<\sigma}{}2$. By the
fact on covering numbers, pick $\sigma<\sigma^{+}<\k$ such that
$\cov(\mu,\k^+,\k^+,\k)=\mu$, and apply theorem (*).

As for no $\mu$ with $\cf\mu=\aleph_0$ is there a $\sigma<\mu$ with
$2^\sigma=\mu$, we can weeken the assumptions to get

\hz{\bf 4.7 Corollary} If $\aleph_{\mu}>\mu$ or $\mu$ is not a second
order fix point, $\cf\mu=\aleph_0$ or $2^{<\cf\mu}<\mu$, and
$\mu\not=2^{<\mu}$, {\it then} there is no universal linear order at
$\mu$.\endproof

\bb
\s{\bf 5 Generalizations}

In this section we prove that if there is no universal linear  order
in a cardinal $\l$ then there is no universal model in $\l$ for any countable
theory $T$ possessing the strict order property (e.g. there is no
universal Boolean algebra in $\l$). This means that all the
non-existence theorems in Section 3 and Section 4 hold for a large
class of theories.
\lz{\bf 5.1 Definition} 1. A formula $\varphi(\ov x;\ov y)$ has the
{\it strict order property} if for every $n$ there are ${\ov a}_l$ ($l<n$)
such that for any $k,l<n$, 
$$
\sat (\exists x)[\neg\phi(\ov x;{\ov a}_k)\wedge \varphi(\ov x; {\ov
a}_l)]\Leftrightarrow k<l
$$
\new 2. A theory $T$ has the {\it strict order property} if some
formula $\varphi(\ov x;\ov y)$ has the strict order property.

This definition appears in in [Sh-a] p. 68, or [Sh-c] p. 69. Every
unstable theory posseses the strict order property or posseses the
independence property (or both). For details see [Sh-a] or [Sh-c].
\lz
{\bf 5.2 Fact} 1. Suppose $T$ has the strict order property, with $\phi(\ov
x;\ov y)$ witnessing this, and let $M\sat T$. Then $\phi$ defines a
partial order $P_M=\lng |M|^n,\le_\phi\rng$, where $n={\rm lg}(\ov y)$, the
length of $\ov y$, the order being given by $\ov y_1\le_\phi \ov y_2
\Leftrightarrow M\sat \phi(\ov x;\ov y_1) \to \phi(\ov x;\ov y_2)$. In
this order there are arbitrarily long chains. 
\hz 2. If $h:M_1\rightarrow M_2$ is an embedding between two models of
$T$ which preserves $\phi$, then $h':P_{M_1}\to P_{M_2}$ is an
embedding of partial orders, where
$h'(x_1,\cdots,x_n)=(h(x_1),\cdots,h(x_n))$.

\new{\ku Proof}: Immediate from the definition.\endproof

\lz {\bf 5.3 Lemma} Suppose $T$ has the strict order property, with
$\varphi(\ov x;\ov y)$ witnessing this. If $L$ is a given linear
order, then there is a model $M$ of $T$ such that $L$ is isomorphic to
a suborder of $P_M$ and $||M||=|L|+\aleph_0$.
\new
{\ku Proof}: As there are arbitrarily long chains with respect to
$\le_\varphi$, this is an immediate corollary of the compactness
theorem. \endproof

\lz
{\bf 5.4 Lemma} If there is a partial order of size $\l$ with the
property that every linear order of size $\l$ can be embedded into it,
then there is a universal linear order in power $\l$.
\hz
{\ku Proof}: Suppose that $P=\lng |P|,\le\rng$ is a partial order of
size $\l$ with this property. Divide by the equivalence relation
$x\sim y \Leftrightarrow x\le y \wedge y\le x$, to obtain $P'$, a
partial order in the strict sense. There is some linear order $<$ on
$|P'|$ which extends $\le$. Let $UO=\lng |P'|,<\rng$. Let $L$ be any
linear order, and let $h:L\to P$ be an embedding. Whenever $x\not= y$
are elements of $L$, $h(x)\not\sim h(y)$ in $P$. Therefore $h':L\to
P'$ defined by $h'(x)=[h(x)]$ is still an embedding. $h'$ is
an embedding of $L$ into $UO$. So $UO$ is universal.\endproof
  
\lz {\bf 5.5 Theorem} Suppose that $T$ has the strict order property,
$\l$ is a cardinal, and that $T$ has a universal model (with respect
to elementary embeddings) in cardinality
$\l$.  Then there exists a universal linear order in cardinality
$\l$.
\hz{\ku Proof}: By 5.4 it is enough to show that there is a
partial order of size $\l$ which is universal for linear orders,
namely that every linear order of the same size is embedded into it.
Let $M$ be a universal model of $T$, and let $\varphi$ witness the
strict order property. We check that $P_M$ is a partial order
universal for linear orders. Suppose that $L$ is some linear order of
size $\l$. By 5.3 $L$ is isomorphic so some suborder of $P_{M_L}$.
Pick an elementary embedding $h:M_L\to M$. In particular, $h$
preserves $\phi$. So by 5.2, there is an embedding of $P_{M_L}$ into
$P_M$. The restriction of this embedding to (the isomorphic copy of)
$L$ is the required embedding.\endproof

\lz{\bf 5.6 Remark} If there is a quantifier free fomula in $T$ which
defines a partial order on models of $T$ with arbitrarily long
chanins, then also a universal model of $T$ in power $\l$ in the sense
of {\it ordinary embedding} implies the existence of a universal
linear order in $\l$.

\lz {\bf 5.7 Conclusions} 
Under the hypotheses of 3.6, there are no universal models in $\l$ for
the following theories:
\hz
Partial orders (ordinary embeddings)
\hz Boolean algebras (ordinary embeddings). 
\hz lattices (ordinary embedding)
\hz ordered fields (ordinary embeddings)
\hz ordered groups (ordinary embeddings)
\hz number theory (elementary embeddings)
\hz the theory of $p$-adic rings (elementary embeddings)

\lz{\ku Proof} All these theories have the strict order property, and
most have a definable order via a quantifier free formula. \endproof

\lz{\bf Appendix}

We review here several motions and definitions from set theory and
model theory.

\new
{\bf Set theory}
\new
A set of ordinals $C$ is {\it closed} if $\sup(C\cap \a)=\a$ implies
$\a\in C$. A {\it club} of a cardinal $\l$ is a closed unbounded set
of $\l$. When $\l$ is an uncountable cardinal, the intersection of two
clubs contains a club. For details see any standard textbook like [Le].

\new
{\bf Model theory}
\new
A model $M$ is a {\it submodel} of $N$ if $|M|\su |N|$ (the universe
of $M$ is a subset of the universe of $N$) and every relation or
function of $M$ is the restriction of the respective relation or
function of $N$. A model $M$ is an {\it elementary submodel} of $N$ if
it is a submodel of $N$, and for every formula $\varphi$ with
parameters from $M$, $M\sat \phi \Leftrightarrow N\sat \phi$. An
embedding of models $h:M\to N$ is an {\it elementary embedding} if its
image is an elementary submodel.  A model $M$ is {\it $\l$-saturated}
if for every set $A\su |M|$ with $|A|<\l$ and a type $p$ over $A$,
$p$ is realized in $M$. A model $M$ is {\it saturated} if it is
$||M||$-saturated. A model $M$ is {\it special} if there is an
increasing sequence of elementary submodels $\lng M_\l:\l<||M||$ is a
cardinal $\rng$, $M=\cup_\l M_\l$ and every $M_\l$ is
$\l^+$-saturated.

A formula $\phi(\ov x;\ov y)$ has the {\it independence property} if
for every $n<\om$ there are sequences $\ov a_l$ $(l<n)$ such that for
every $w\su n$, 
$\sat (\exists\ov x)\wedge_{l<n}\phi(\ov x;\ov a_l)^{{\rm if}l\in w}$
\new
A first order theory $T$ has the {\it independence property} if some
formula $\phi(\ov x;\ov y)$ has the independence property.

\new{\bf A theory with universal models with respect to ordinary but
not elementary embeddings}

Let $T$ be the theory of the following model $M$: the universe is
divided into two infinite parts, the domain of the unary predicate $P$
and its complement. The domain of $P$ is the set $\lng
x_i:i<\om+\om\rng$. There are two binary relations $R_1, R_2$. For
every pair $\lng \a,\b\rng$ of ordinals below $\om+\om$ there is a
unique element $y\in M$ such that $\neg P(y)$ (think of $y$ as an ordered
paire) which satisfies $R_1(y,x_\a)\wedge R_2(y,x_\b)$ iff $\a<\b<\om$
or $\a<\om\le \b$ or $\a\ge \om$ and $\b\ge \om$. So for the elements
above $\om$ in the $P$ part all possible ordered pairs exist, while
those below $\om$ are linearly ordered by the existence of ordered
pairs. Let $\phi(x_1;x_2)=P(x_1)\wedge P(x_2)\wedge (\exists y)(\neg
P(y)\wedge R_1(y,x_1)\wedge R_2(y,x_2)$. This formula witnesses that
$T$ has the strict order property. By Theorem 3.5 and Theorem 5.5 $T$
has no universal model with respect to elementary embeddings in any
regular $\l$, $\aleph_1<\l<\aleph_0$. But with respect to ordinary
embeddings $T$ has a universal model in every infinite cardinality
$\l$: let $M$ be any model of $T$ of cardinality $\l$ in which there
is a set $X=\{ x_i:i<\l \}$ of elements in the domain of $P$ for which
all possible ordered pairs exist. If $M'$ is any other model of
cardinality $\l$, and $h$ is a 1-1 function which maps the domain of
$P$ in $M'$ into $X$, then $h$ can be completed to an
embedding of $M'$ into $M$.\endproof

\new{\bf The consistency of not having universal models}
 Let us state it here, for the sake of those who
read several times that was is easy but are still interested in the
details:
\hz
{\bf Fact} If $\l$ regular, $V\sat$ GCH, for simplicity, and $P$ is
a Cohen forcing which adds $\l^{++}$ Cohen subsets to $\l$, then in $V^P$
there is no universal graph (linear order, model of a complete first order
 $T$ which is unstable in $\l$) in power $\l^+$.
\hz{\ku Proof}
Suppose to the contrary that there is a universal graph $G^*$ of power
$\l^+$. We may assume, without loss of generality, that its universe
$|G^*|$ equals $\l^+$. As $G$ is an object of size $\l^+$, it is in
some intermediate universe $V'$, $V\su V'\su V^P$, such that there are
$\l^{++}$ Cohen subsets outside of $V'$. So without loss of generality,
$G\in V$. Let $G$ be the graph with universe $\l^+$
defined as follows: fix a 1-1 enumeration $\lng A_\a:\l\le\a<\l+\rng$ of
$\l^+$ Cohen subsets of $\l$. A pair $\a,\b$ is joined by an edge iff
$\b\ge\l$ and $\a\in A_\b$,  or $\a\ge\l$ and $\b\in A_\a$. By the
universality of $G^*$, there should exist an embedding $h:G\rightarrow
G^*$. Consider $h\rest\l$. This is an object of size $\l$, and
therefore is is some $V'$, $V\su V'\su V^P$, where in $V'$ there are
at most $\l$ of the Cohen subsets. For every $y\in G^*$, the set
$\{x\in\l:h(x)$ is joined by an edge to $y\}$ is in $V'$. Pick an
$\a\ge\l$ such that $A_\a$ is not in $V'$ and set $y=h(\a)$ to get a
contradiction.\endproof

The same proof is adaptable to the other cases.

\new {Guessing clubs}
\hz{\ku Proof of Fact 3.3}
Let $S_0$ be $\{\d<\l:\cf\d=\aleph_0\}$. Suppose to the contrary that
for every sequence $\ov C$ as above there is a club $C\su\l$ such that
for every $\d\in S_0\cap C$, $c_\d\not\su C$. We construct by
induction on $\b<\al_1$ $\ov C_\b, C_\b$.  $\ov C_\b=\lng
c^\b_\d:\d\in S_0\rng$ is such that for every $\d\in C_\b\cap S_0$,
$c^\b_\d$ is a club of $\d$ of order type $\om$, and $c^\b_\d\not\su
C_\b$.  Furthermore, letting $c^\b_\d=\lng
\a^{\b,\d}_n:n<\om\rng$, where $\a^{\b,\d}_n\le\a^{\b,\d}_m$ if $n<m$,
we demand that $\a^{\b+1,\d}_n=sup\{\a^{\b,\d}_n\cap C_\b\}$ if this
intersection is non-empty, and $\a^{\b+1,\d}_n=0$ otherwise. When $\b$
is limit, we demand that $\a^{\b,\d}_n=min\{\a^{\c,\d}_n$ for
$\c<\b\}$.  Let $\ov C_0$ be arbitrary. At the induction step pick
$C_\b$ which demonstrate that $\ov C_\b$ is not as required by the
fact and define each $c^{\b+1}_d$ by the demand above. Note that for
club many $\d$'s the resulting $c^{\b+1}_\d$ is cofinal in $\d$, so
without loss of generality this is so for every $\d\in C'_{\b+1}$.

It is straightforwaed to verify that for all $\b<\c<\om_1$, $\d\in S_0$
\item {1.}  For  all $\d\in S_0\cap C_\b$, $c^{\b,\d}\su \d$ is a club of
$\d$ of order type $\om$.
\item {2.}  For all $\d\in S_0$, $c^\b_\d\sm\{0\}\su C_{\b+1}$.
\item {3.}  For all $\d\in S_0\cap C_\b$, $c^\b_\d\not \su C_\b$
\item {4.}  $\a^{\c,\d}_n\le\a^{\b,\d}_n$

Let  $C=\bigcap_{\b<\om_1}C_\b$. Pick
$\d_0\in C\cap S_0$. Then $C\cap \d$ is unbounded in $\d$ and
of order type $\om$. Furthermore, for every $\b<\om_1$,
$c^\b_{\d_0}\not\su C_\b$. 
\vskip.1cm
But on the other hand, there is a $\b_0$
such that for all $\b_0<\b<\c$, $c^\b_{\d_0}=c^\c_{\d_0}$, because of
4. - a contradiction.\endproof

Now for the case of uncountable cofinality. We recall 
\hz 
{\bf 3.4 Fact}: If $\cf\k=\k<\k^+<\l=\cf\l$, then there is a
sequence $\ov C=\lng c_\d:\d\in S^\l_\k\rng$, $S^\l_\k$ being the set
of members of $\l$ with cofinality $\k$, where $c_\d$ is a club of $\d$
of order type $\k$ with the property that for every club $E\su \l$ the
set $S_E=\{\d\in S^\l_\k:c_\d\su E\}$ is stationary. Furthermore, if there
is a square sequence on $S^\l_\k$, $\ov C$ can be chosen to be a
square sequence.
\hz {\ku Proof}: This proof is actually simpler than that of the
previous fact. Start with any sequence $\ov C_0=\lng c^0_\d:\d\in
S^\l_\k\rng$. By induction on $i<\k^+$ define $\ov C_i$ as follows: if
$\ov C_i$ has the property of guessing clubs, we stop. Otherwise there
is a club $E_i$ such that $E_i$ is not guessed stationarily often.
This means there is a club $C_i$ such that $\d\in C_i$ implies that
$c_{i_\d}\not\su E_i$. We may assume that $C_i=Ei$. Let
$c^{i+1}_\d=c^i_\d\cap E_i$. If $i$ is limit,
$c_\l^i=\cap_{j<i}c^j_\d$. Suppose the inductions goes on for $\k^+$
steps. Let $E=\cap E_i$. For every $\d\in E\cap S^\l_\k$, $E\cap
C_\d^0$ is a club of $\d$. Therefore for stationarily many points
$\d$, $c^{\k^+}$ is a club of $\d$, say that this holds for all $\d\in
S\su S^\l_\k$. Clearly, for every $\d$ there is an $i$ such that for
every $i<j$ $C^i_\d=c^j_\d$, because the size of $c^0_\d$ is $\k$ and
$\k^+$ is regular. But on the other hand, for every $\d\in S$ and
$\i<\k^+$, $c^i_\d\not \su E_i$, while $c^{i+1}_\d\su E_i$. A
contradiction. Thus the induction stops before $\k^+$, and the
resultiong sequence guesses clubs.

 What about the square property?
If there is a square sequence on $S^\l_\k$, let in the proof $\ov C_0$
be a square sequence. Notice that the operation of intersecting the $c_\d$ with a club $E$
preserves the property: suppose $\d_1<\d_2$ amd
$\d_1\in\acc(\c_{d_2})\cap E$. The clearly $\d_1\in \acc (c_{\d_2})$.
Therefore, by the square property, $c_{\d_1}=\c_{\d_2}\cap \d_1$.
Intersecting both sides of the equation with $E$ yields that
$c_{\d_1}\cap E=c_{\d_2}\cap E \cap \d_1$. Therefore the proof is
complete. \endproof
\bigskip
\centerline{\bf References}

\lz
[Le] Azriel Levy, {\bf Basic Set Theory}, Springer Verlag, 1979.
\lz
[CK] C. C Chang and J. J. Keisler {\bf Model Theory}, North Holland 1973. 
\lz
[Sh-a] Saharon Shelah, {\bf    Classification Theory and the Number of
non-Isomorphic Models}, North Holland, 1978.
\lz
[Sh-c] Saharon Shelah, {\bf Classification Theory and the Number of
non-Isomorphic Models, Revised}, North Holland Publ. Co., 1990.
\lz
[GrSh 174] R. Grossberg and S. Shelah, {\sl On universal locally
finite groups}, Israel J. of Math. {\bf 44 }, (1983), 289-302.
\lz
[Sh-g] Saharon Shelah, {\bf Cardinal Arithmetic}, to appear.
\lz
[Sh 93] Saharon Shelah, {\sl Simple Unstable Theories}, Annals of Math.
Logic {\bf 19} (1980), 177--204.
\lz
[Sh 175] Saharon Shelah, {\sl On Universal graphs without instances of CH} Annals of
Pure and Applied Logic {\bf 26},
(1984), 75--87.
\lz
   [Sh 175a] Saharon Shelah, {\sl Universal graphs without instances of CH, revised},
Israel J. Math, {\bf 70} (1990) 69--81.
\lz
[M] Alan H. Mekler, {\sl Universal Structures in Power $\aleph_1$}
Journal of Symbolic Logic, accepted.

  \end